\documentclass[a4paper, 12pt]{amsproc}
\usepackage{amssymb,amsthm,amsfonts,latexsym}
\newtheorem{theorem}{Theorem}[section]

\theoremstyle{remark}

\title[Canonical Left Cells and the Lowest Two-sided Cell]
{Canonical Left Cells and the Lowest Two-sided Cell in an Affine
Weyl Group}

\dedicatory{Dedicated to Professor George Lusztig on his seventieth
birthday}
\author[N. Xi]{ Nanhua Xi$^{*}$}
\address{$^{*}$
Institute of Mathematics\\
Chinese Academy of Sciences\\
Beijing, 100190\\ and\\ School of Mathematical Sciences, University of Chinese Academy of Sciences\\ Beijing 100049\\
China } \email{nanhua@math.ac.cn}
\thanks{This work is partially supported by Natural Sciences Foundation
of China (No. 11321101).}

\begin{document}
\baselineskip=18pt
\begin{abstract}
We give
some discussions to the relations between  canonical left cells and
the lowest two-sided cell of an affine Weyl group. In particular, we
use the relations to construct irreducible modules attached to the lowest two-sided cell and some one dimensional representations
of an affine Hecke algebra.

\end{abstract}

\maketitle

\def\Cal{\mathcal}
\def\bold{\mathbf}
\def\ca{\mathcal A}
\def\cdz{\mathcal D_0}
\def\cd{\mathcal D}
\def\cdo{\mathcal D_1}
\def\bold{\mathbf}
\def\l{\lambda}
\def\le{\leq}
\def\bbZ{\mathbb Z}

Canonical left cells of an affine Weyl group are interesting in
understanding cells in affine Weyl group and have nice relations
with structure and representations of algebraic groups. However, it
is not easy to describe canonical left cells. In this paper we give
some discussions to the relations between  canonical left cells and
the lowest two-sided cell of an affine Weyl group. In particular, we
use the relations to construct irreducible modules attached to the lowest two-sided cell (see Theorem4.1)and some one dimensional representations
of an affine Hecke algebra (see Theorem 3.5). For convenience we work with an extend
affine Weyl group. This work was partially motivated by [AB].

\section{Canonical left cells}

\subsection{}
  Let $R$ be an irreducible root system and $P$ the corresponding
  weight lattice.
  The Weyl group $W_0$ acts on $P$ naturally and the semi-direct product
$W=W_0\ltimes P$ is an extended affine Weyl group, which contains
the affine Weyl group $W_a=W_0\ltimes \bbZ R$. Let $S$ be the set of
simple reflections of $W_a$. The partial order $\le $ and the length
function $l$ on $W$ are well defined. {\it The operation on $P$ will be written in multiplication.}

For $w\in W$, set  $L(w)=\{s\in S\,|\, sw\le w\}$ and $R(w)=\{s\in
S\,|\, ws\le w\}$.
Let $s_0$ be the unique simple reflection of $W_a$ out of $W_0$.
Define $Y_0=\{w\in W\,|\, R(w)\subseteq\{s_0\}\,\}$.  Then
  $Y_0\cap \Omega$ is a left cell for any two-sided cell $\Omega$ of
  $W$, called a canonical left cell.

In general it is not easy to describe a canonical left cell.
However, it is  easy to describe the set $Y_0$. Let $w_0$ be the
longest element of $W_0$. The set of anti-dominant weights in $P$ is
defined to be  $P^-=\{x\in P\,|\, l(xw_0)=l(w_0)+l(x)\}$ and the set
of dominant weights is $P^+=\{x\in P\,|\, l(w_0x)=l(w_0)+l(x)\}$.

\smallskip

\noindent {\bf Proposition 1.2.} $Y_0=\{wx\,|\, w\in W_0,\ x\in P^-
\text{ and } R(w)\subseteq L(x)\}$.

\smallskip

Proof. Let $u\in W$. Then there exist unique $w,v\in W_0$ and $x\in
P^-$ such that $R(w)\subseteq L(x)$ and $u=wxv$. Moreover, we have
$l(u)=l(x)+l(v)-l(w)$. The proposition follows.

\smallskip

\noindent{\bf 1.3.} It would be interesting to see when two elements
in $Y_0$ are in a left cell. Let $\rho$ be the product of all
fundamental dominant weights. Then the set $\{wx\rho^{-1}\,|\, w\in W_0,\
x\in P^- \text{ and } R(w)\subseteq L(x)\}$ is the canonical left
cell in the lowest two-sided cell $c_0$ of $W$. In general, for any
$x\in P^-$ there exists a positive integer $a$ (depending on $x$)
such that $x^b$ and $x^a$ are in a left cell if $b\ge a$ (see [X1,
Lemma 3.2]). It seems that the number $a$ is not big, in many cases,
it is among 1,2,3.

\def\Gz{\Gamma_0}

Let $S_0=S\cap W_0$ and denote by $\Gamma_0$ the left cell $\{w\in
W\,|\, R(w)=S_0\}$, which is in the lowest two-sided cell $c_0$ of
$W$.  For $x\in P$, denote by  $n_x$ (resp. $m_x$)  the unique
shortest element in the coset $xW_0$ (resp. the double coset
$W_0xW_0$). The map $x\to n_x$ defines a one-to-one correspondence
from $P$ to $Y_0$, and the map $n_x \to n_xw_0$ defines a one-to-one
correspondence from $Y_0$ to $\Gz$. Also the map $x\to m_x$ defines
a one-to-one correspondence between $P^+$ and $Y_0\cap Y_0^{-1}$.
 The sets $ Y_0$ and
$\Gz$ produce naturally two modules of an affine Hecke algebra of
$(W,S)$. In next section we will see that the two modules are
essentially the same.

\def\a{\alpha}

\section{Cell modules of affine Hecke algebras}

\noindent{\bf 2.1.}  Let $H$ be the  Hecke algebra  of $(W,S)$ over
a field $k$ with parameter $q$. Assume that $k$ contains
 square roots of $q$.  Let  $\{T_w\}_{w\in W}$ be its standard basis.
For any $w$ in $W$, let $$C_w=q^{-\frac
{l(w)}2}\sum_{y\le w}P_{y,w}(q)T_y,$$ and $$C'_w=q^{\frac {l(w)}2}\sum_{y\le
w}(-1)^{l(w)-l(y)}q^{-l(y)}P_{y,w}(q^{-1})T_y,$$ where $P_{y,w}$ are the
Kazhdan-Lusztig polynomials. Then the elements $C_w,\ w\in W$, form a
basis of $H$, and the elements $C'_w, w\in W$, form a  basis of $H$
as well, see [KL1].

\def\tx{\theta_x}
\def\ty{\theta_y}
\def\t{\theta}
\def\T{\Theta}

For any $x\in P$ there is a well defined element
$\theta_x=q^{-\frac{l(y)}2}T_y q^{\frac{l(z)}2}T_z^{-1}$. where
$y,z\in P^+$ such that  $x=yz^{-1}$.  Then $\theta_x\theta_y=\ty\tx$
for any $x,y\in P$ and the elements $T_w\tx$ (resp. $\tx T_w$), $ w\in
W_0,$ $x\in P$, form a basis of $H$. See [L1].

The group algebra $k[P]$ is isomorphic to the subalgebra $\Theta$ of
$H$ generated by all $\tx,\ x\in P$. Lusztig defined several
$H$-module structures on $k[P]$, see [L2, Section 7]. They are
actually isomorphic to the modules provided by the left cell $\Gz$. Let
$M$ (resp. $M'$) be the subspace of $H$ spanned by all $C_w,\
w\in\Gz$ (resp. $C'_w,\ w\in\Gz$). Then $M$ and $M'$ are left
ideals of $H$ and generated by $C=C_{w_0}$ and $C'=C'_{w_0}$
respectively. The elements $\tx C$, $x\in P$, form a basis of $M$
and the elements $\tx C'$, $x\in P$, form a basis of $M'$.

Let $\mathcal I$ (resp. $\mathcal I'$) be the subspace of $H$ spanned by all $C_w, \
w\in W-Y_0$ (resp. $C'_w,\ w\in W-Y_0$). Then $\mathcal I$ and $\mathcal I'$ are left
ideals of $H$. Let $N=H/\mathcal I$ and $N'=H/\mathcal I'$. Essentially the following
result is due to Arkhipov and Bezrukavnikov (see [AB, 1.1.1]).

\medskip

\noindent {\bf Lemma 2.2.}  As $H$-modules $N$ is isomorphic to
$M'$, and $N'$ is isomorphic to $M$.

\medskip

Proof. Consider the surjective homomorphism $H\to M'$, $h\to hC'$.
It is easy to check that the kernel is $\mathcal I$. So $N$ is isomorphic to
$M'$. Similarly the surjective homomorphism $H\to M$, $h\to h C$
induces an isomorphism $N'\to M$ of $H$-module. The lemma is proved.

\medskip

\noindent{\bf 2.3.} The geometric explanation of the isomorphism in the above lemma is
that Thom isomorphism for a certain equivariant $K$-group  of the
cotangent bundle of flag variety is compatible with certain actions
of the affine Hecke $H$, see [L2, Section 7].

Lemma 2.2 seems helpful in understanding the structure of
$H$-modules $M$ and $M'$, and may be useful to understand canonical
left cells. A natural question is to consider the submodule of $M'$
(resp. $M$) generated by  all $C_wC'$ (resp. $C'_wC$),  $w\in
c_0\cap Y_0$. Modulo a central character of $H$, we can get a finite
dimensional quotient algebra of $H$. In next two sections we will give
some discussion to the images in such quotient algebras of the
submodules. We will show that the images in such a quotient algebra is  either irreducible $H$-module  or zero when $k$ is algebraically closed (Theorem 4.1).

\section{A realization of some one dimensional representations}

In this section we construct some one dimensional representations of the affine Hecke algebra $H$ through certain quotient algebras of $H$ (see Theorem 3.5).

\subsection{}  From now on we assume that $k$ is algebraically closed.  Recall that $\Theta$ is the subalgebra of $H$ generated by all $\theta_x,\ x\in P$. Let the Weyl group $W_0$ act on $\Theta$ by $w(\theta_x)=\theta_{w(x)}$.

 We shall need a several formulas in $H$. Let $x\in P$, the Macdonald formula says (see [NR,
Theorem 2.22])
\begin{equation} C\tx C=q^{-\frac{l(w_0)}2}\sum_{w\in W_0}w(\tx \prod_{\alpha\in
R^+}\frac{1-q\t_{\alpha}}{1-\t_{\alpha}})C.\end{equation}

 Let $\Delta$ be the set of simple roots of $R$ and
 denote $x_\a$ the fundamental dominant weight corresponding to a simple root $\alpha$.
Recall that $\rho=x_{\Delta}$ is the product of all fundamental dominant
weights. Using [L1, Corollary 7.8, Lemma 7.4 (iii)] we get
\begin{equation}C'\t_{\rho^{-1}} C= q^{-\frac{\nu}2}C' \sum_{I\subseteq
R^+}(-q)^{|I|}\t_{\rho^{-1}}\t_{\a_I},\end{equation}
\begin{equation}C'\t_{\rho} C= q^{\frac{\nu}2}C' \sum_{I\subseteq
R^+}(-q)^{-|I|}\t_{\rho}\t_{\a_I^{-1}},\end{equation}where $\nu=l(w_0)=|R^+|$,
$\a_I$ is the sum of all roots in $I$ and $|I|$ is the cardinality
of $I$.

There is a unique involutive anti-automorphism $h\to \tilde h$ of
the $k$-algebra $H$ such that $\tilde T_r=T_r\ (r\in S_0), \
\tilde\tx=\t_{x}\ (x\in P)$ [KL2, 2.13(c)].  Note that $\tilde C=C'$
and $\tilde C'=C$. Applying this anti-automorphism to the formulas (2) and (3) we get
\begin{equation}C\t_{\rho^{-1}} C'= q^{-\frac{\nu}2} \sum_{I\subseteq
R^+}(-q)^{|I|}\t_{\rho^{-1}}\t_{\a_I}C',\end{equation}
\begin{equation}C\t_{\rho} C'= q^{\frac{\nu}2} \sum_{I\subseteq
R^+}(-q)^{-|I|}\t_{\rho}\t_{\a_I^{-1}}C'.\end{equation}
(For a K-theoretic understanding of formula (4) see [X6, 2.6], note that the $C'$ here is the $C$ in loc.cit.)

\subsection{The center of $H$}  We have

\medskip

(a) The center $Z(H)$ of $H$  consists of $W_0$-invariant elements in $\Theta$, i.e. $Z(H)=\Theta^{W_0}$ (see [L1]).

\medskip

Therefore the
center $Z(H)$ of $H$  is isomorphic to $k\otimes_{\bbZ} R_G$, where $G$ is a simply connected
simple algebraic group over $k$ with root system $R$ and $R_G$ is
the representation ring of $G$.

\def\st{\stackrel}
\def\sc{\scriptstyle}
For $w\in W_0$ define
$$e_w=w(\prod_{\st {\alpha\in\Delta}{ w(\alpha)\in R^-}}x_\alpha).$$ (Recall that here $x_\a$ is the fundamental dominant weight
corresponding to $\alpha\in\Delta$.) Then (see [S])

\medskip

(b) $\T$ is a free
$Z(H)$-module with a basis $\{\t_{e_w}\,|\, w\in W_0\}$.

 \medskip

 Hence,

 \medskip

 (c) for any $x\in P$, the elements $\t_x\t_{e_w},\ w\in W_0,$ form a $Z(H)$-basis of $\T$.

\medskip

(d) For $A,B\in \T$, define
$$(A,B)=(-1)^\nu\t_\rho\prod_{\a\in R^+}(1-\t_\a)^{-1}\sum_{w\in
W_0}(-1)^{l(w)}{w(AB\t_\rho)}\in Z(H).$$ By [KL2, p.163] there exist
$\t'_u\in\T\  (u\in W_0)$ such that $(\t_{e_w},\t'_u)=\delta_{w,u}$
and the elements $\t'_u$ form a $Z(H)$-basis of $\T$.

\subsection{}  The set of $k$-algebra
homomorphisms from $Z(H)$ to $k$ is  in one-to-one correspondence to
the set of semisimple classes of $G$. For each semisimple class
$\bar t$ in $G$, let $\phi_{\bar t}: Z(H)\to k$ be the corresponding
homomorphism. Let $T$ be a maximal torus of $G$ and identify $P$
with the character group Hom$(T,k^*)$ of $T$. For $t\in T$, we have
a $k$-algebra homomorphism $\phi_t:\T\to k$ defined by $\tx\to x(t)$
for all $x\in P$. Let $\bar t$ be the conjugacy class of $G$
containing $t$. Then the restriction to $Z(H)$ of $\phi_t$ is
$\phi_{\bar t}$.   For a semisimple element $t$ in $T$, let $\mathcal Z_t$ be
the two-sided ideal of $H$ generated by all $z-\phi_t(z),\ z\in
Z(H)$. Define $H_t=H/\mathcal Z_t$. Then dim$H_t=|W_0|^2$.

For each simple $H$-module $L$, there exist some $t$ in $T$ such
that $Z(H)$ acts on $L$ through the homomorphism $\phi_t$. So to
study simple modules of $H$ it is enough to study simple modules of
the quotient algebras $H_t$ for $t\in T$. We shall use the same
notations $C'_w, D_w, C, C', \tx, ...$  for their images in $H_t$.

\renewcommand{\thetheorem}{3.4}
\begin{theorem}Let $t\in T$. The following statements
are equivalent.

(a) $CH_{t}C=0.$ (Recall that $C=C_{w_0}$ and $C'=C'_{w_0}$.)

(b) $CH_tC'=0.$

(c)  $C'H_tC=0.$

(d)  $C'H_{t^{-1}}C'=0.$

(e) For any simple $H_{t}$-module $L$ we have $CL=0$.

(f) For any simple $H_{t^{-1}}$-module $L$ we have $C'L=0$.\end{theorem}

\smallskip

Proof.  There is a unique involutive automorphism $h\to h^*$ of the
$k$-algebra $H$ such that $T_r^*=-qT_r^{-1}=q-1-T_r\ (r\in S_0), \
\tx^*=\t_{x^{-1}}\ (x\in P)$ [KL2, 2.13(d)]. Noting that
$C^*=(-1)^{l(w_0)}C'$, we see that (a) and (d) are equivalent, (e)
and (f) are equivalent.

Using the  involutive anti-automorphism $h\to \tilde h$ of
the $k$-algebra $H$ defined by $\tilde T_r=T_r\ (r\in S_0), \
\tilde\tx=\t_{x}\ (x\in P)$ [KL2, 2.13(c)] and noting that $\tilde C=C$
and $\tilde C'=C'$, we see that (b) and (c) are equivalent.

Since the two-sided ideal $H_{c_0}$ of $H$ spanned by all $C_w,\
w\in c_0$ is generated by $C$, using [X4, 7.7] we know that (a) and
(e) are equivalent.

Now we show that (a) and (b) are equivalent. Since
$T_wC=CT_w=q^{l(w)}C$ if $w\in W_0$, we see that $CHC$ is spanned
by $C\tx C$. By formula (1) in 3.1, we have
$$C\tx C=q^{-\frac{l(w_0)}2}\sum_{w\in W_0}w(\tx \prod_{\alpha\in
R^+}\frac{1-q\t_{\alpha}}{1-\t_{\alpha}})C.$$ So we have

\noindent (i) The condition $CH_tC=0$ is equivalent to
 $$\phi_t(\sum_{w\in W_0}w(\tx \prod_{\alpha\in
R^+}\frac{1-q\t_{\alpha}}{1-\t_{\alpha}}))=0,\quad\text{for all
}x\in P. $$

\def\a{\alpha}

\medskip

Using 3.2 (c) we know that $H$ is spanned by all $T_wz\t_\rho\t_{e_u},$ $w,u\in W_0,\ z\in Z(H)$. Therefore we have

\noindent(ii) The condition $CH_tC=0$ is equivalent to
 $$\phi_t(\sum_{w\in W_0}w(\t_\rho\t_{e_u} \prod_{\alpha\in
R^+}\frac{1-q\t_{\alpha}}{1-\t_{\alpha}}))=0,\quad\text{for all
}u\in W_0. $$

 Similar to [X3,  Lemma 2.10], we see that $HC'$ is
spanned by all $T_wz\t_{I}C'$, $w\in W_0,\ z\in Z(H),\ I\subseteq
\Delta$, where $\t_I=\prod_{\a\in I}\t_{x_\a}$.
Since $T_wC=CT_w=q^{l(w)}C$  if $w\in W_0$ and $C\t_IC'=0$ if $I\ne
\Delta$, as a $Z(H)$-module, $CHC'$ is generated by $C\t_\rho C'$.

Recall the formula (5) in 3.1:
$$C\t_\rho C'= q^{\frac{\nu}2}\sum_{I\subseteq
R^+}(-q)^{-|I|}\t_\rho\t_{\a_I^{-1}}C',$$ where $\nu=l(w_0)=|R^+|$,
$\a_I$ is the sum of all roots in $I$ and $|I|$ is the cardinality
of $I$.

\def\st{\stackrel}
\def\sc{\scriptstyle}

Let $A=q^{\frac{\nu}2}\sum_{I\subseteq R^+}(-q)^{-|I|}\t_\rho\t_{\a_I^{-1}}$. Note
that $$A=(-1)^\nu q^{-\frac{\nu}2}\t_\rho^{-1}\prod_{\a\in R^+}(1-q\t_\a).$$ Thus
\begin{equation}(A,\t_{e_u})=(-1)^\nu q^{-\frac{\nu}2}\sum_{w\in W_0}w(\t_{\rho e_u}\prod_{\a\in
R^+}\frac{1-q\t_\a}{1-\t_\a}).\end{equation} Since $A=\sum_{u\in
W_0}(A,\t_{e_u})\t'_u$ in $H$ and $\theta'_uC',\ u\in W_0$, are linearly independent in $H_tC'$, we obtain

\noindent(iii) The condition $CH_tC'=0$ is equivalent to
$$\phi_t(\sum_{w\in W_0}w(\t_{\rho e_u}\prod_{\a\in
R^+}\frac{1-q\t_\a}{1-\t_\a})=0\quad\text{for all }u\in W_0.$$

Using (ii) and (iii) we see that (a) and  (b) are equivalent.  The theorem is proved.

\renewcommand{\thetheorem}{3.5}
\begin{theorem}
 Let $t\in T$ be such that $\alpha(t)=q$
for all simple roots $\alpha$ of $R$. Then

(a) $CH_tC'$ (resp. $C'H_tC$) is a two-sided ideal of $H_t$ with
dimension 1 if $\sum_{w\in W_0}q^{l(w)}\ne0.$

(b) $CH_tC'=0$ if $\sum_{w\in W_0}q^{l(w)}=0.$\end{theorem}

\smallskip

Proof. We have seen that $CH_tC'$ is spanned by the image in $H_t$
of $C\t_\rho C'$. To see it is a two-sided ideal of $H_t$ it
suffices to prove that the images in $H_t$ of $\tx C\t_\rho C'$ and
$ C\t_\rho C' \tx$ for all $x\in \T$ are scalar multiples of the
image in $H_t$ of $C\t_\rho C'$.

\noindent (i) If $w$ is not the neutral element of $W_0$, then there
exists a positive root $\beta$ such that $w(\beta)=\a^{-1}$ for some
simple root $\a$. Thus $w(1-q\beta)(t)=0$.

Let $A$ be as in the proof  of Theorem 3.4. Then $\t_xC\t_\rho C'=A\t_xC'$.
Since
$$(A\tx, \t_{e_u})=(-1)^\nu q^{-\frac{\nu}2}\sum_{w\in W_0}w(\t_{\rho xe_u}\prod_{\a\in
R^+}\frac{1-q\t_\a}{1-\t_\a}),$$ using (i) we get
$$(A\tx,\t_{e_u})(t)=\rho(t)x(t)e_u(t)\prod_{a\in R^+}\frac{1-q^{1+\langle
\rho,\a^\vee\rangle}}{1-q^{\langle \rho,\a^\vee\rangle}},$$ if
$1-q^{\langle \rho,\a^\vee\rangle}\ne 0$ for all positive roots
$\alpha$. We have (see for example [NR, Corollary 2.17])
$$\prod_{a\in R^+}\frac{1-q^{1+\langle
\rho,\a^\vee\rangle}}{1-q^{\langle \rho,\a^\vee\rangle}}=\sum_{w\in
W_0}q^{l(w)},$$ if $1-q^{\langle \rho,\a^\vee\rangle}\ne 0$ for all
positive roots $\alpha$. Now $(A\tx,\t_{e_u})$ is in $Z(H)$, so
$(A\tx,\t_{e_w})(t)$ is a regular function in $q\in k^*$. Thus we
have \begin{equation}(A\tx,\t_{e_u})(t)=(-1)^\nu q^{-\frac{\nu}2}\rho(t)x(t)e_u(t)\sum_{w\in W_0}q^{l(w)}\end{equation} for all
$q\in k^*$. So the images in $H_t$ of $\tx C\t_\rho C'$ for all $x\in
\T$ are scalar multiples of the image in $H_t$ of $C\t_\rho C'$, and
$CH_tC'$ is a left ideal of $H_t$. Using the
involutions $h\to h^*$ and $h\to \tilde h$ of $H$ several times we see that
$CH_tC'$ is a left ideal of $H_t$ implies that it is also a right
ideal of $H_t$.

The formula (7) also indicates that $CH_tC'=0$ if and only if $\sum_{w\in W_0}q^{l(w)}=0$. The theorem is proved.

\smallskip

It is easy to check that $T_sCH_tC'=qH_tC'$ and $CH_tC'T_s=-CH_tC'$
for all simple reflections $s$ if $\a(t)=q$ for all simple roots
$\a$. So the ideals $CH_tC'$ and $C'H_tC$ give natural realizations
of some one dimensional representations of $H_q$.

\section{Irreducible modules attached to the lowest two-sided cell}

The main result of this section is the following.

\renewcommand{\thetheorem}{4.1}
\begin{theorem} Let $t\in T$, then

(a) The element $C\theta_\rho C'$ in $H_t$  generates an irreducible module $L_t$ of $H$ if it is nonzero. Moreover, $CL_t\ne 0$ in this case.

(b) The element $C'\theta_\rho C$ in $H_t$ generates an irreducible module $L'_t$ of $H$ if it is nonzero. Moreover, $C'L_t\ne 0$ in this case. \end{theorem}

Proof. Let $J_{c_0}$ be the based ring of $c_0$. According to [X2], $J_{c_0}$ is isomorphic to a $|W_0|\times |W_0|$  matrix ring over $R_G$. Let $\mathbf J_{c_0}=\mathbb C\otimes J_{c_0}$. Then up to isomorphism,  irreducible $\mathbf J_{c_0}$-modules are naturally one-to-one corresponding to the semisimple classes of $G$. For the semisimple class containing $t$, let $E_t$ be a corresponidng simple $\mathbf J_{c_0}$-module.

Let $\varphi_0:H\to \mathbf J_{c_0}$ be Lusztig's homomorphism defined through the basis $C_w, w\in W$. Then $E_t$ is endowed with an $H$-module structure through the homomorphism. Denote the $H$-module structure on $E_t$ by $E_{t,\varphi_0}$. We have (see for example the proof 2.5 in [X6])

(i) $E_{t,\varphi_0}$ is isomorphic to $H_tC$.

According to [X4, 7.7] and [X5, Lemma 2.5], $H_tC$ has a simple constituent $L$ such that $CL\ne 0$ if and only if $CH_tC\ne 0$. In this case, $L$ is the unique simple simple constituent of $H_tC$ such that $CL\ne 0$ and $L$ is also the unique simple quotient module of $H_tC$.

By Theorem 3.4, $CH_tC\ne 0$ is equivalent to $CH_tC'\ne 0$. Now assume that $CH_tC'\ne 0$. By the proof of [X5, Lemma 2.5], the set
$$M_{t,0}=\{ h \in H_tC\,|\, C_wh=0,\quad \forall w\in c_0\}$$
is the unique maximal submodule of $H_tC$.

We have a  natural $H$-module homomorphism:  $H_tC\to H_tC\t_\rho C'$, $h\to h\t_\rho C'$. Therefore, to prove that $C\theta_\rho C'$ in $H_t$  generates an irreducible module $L_t$ of $H$ it suffices to prove that
$$h\in M_{t,0}\Longleftrightarrow h\t_\rho C'=0.$$

Let $\T_t$ be the image of $\T$ in $H_t$. Then $H_tC$ consists of  $\theta C,\ \t\in\T_t$. Since $\T$ is a free $Z(H)$-module with a basis $\{e_w\,|\, w\in W_0\}$, and for any $w\in c_0$ there exists $\xi,\ \eta\in \T$ such that $C_w=\xi C\eta$, we see that $\t C\in H_t C$ is in $M_{t,0}$ if and only if $C\t_{e_w}\t C=0$ in $H_t$ for all $w\in W_0$. By formula (1) in 3.1, this is equivalent to
\begin{equation}\phi_t(\sum_{w\in W_0}w(\t_{e_u}\t \prod_{\alpha\in
R^+}\frac{1-q\t_{\alpha}}{1-\t_{\alpha}}))=0,\quad\text{for all
}u\in W_0. \end{equation}

Let $A$ be as in the proof of Theorem 3.4, then $$\t C\t_\rho C'=\t AC'.$$
Clearly, $\t C\t_\rho C'=0$ in $H_t$ if and only if  $\t_{\rho}^{-1}\t C\t_\rho C'=0$ in $H_t$.
Since $$\t_{\rho}^{-1}\t C\t_\rho C'=\sum_{u\in W_0}(\t_{\rho}^{-1}\t A,\t_{e_{u}})\t'_uC'$$ and  $\theta'_uC',\ u\in W_0$, are linearly independent in $H_tC'$, we see that $\t C\t_\rho C'=0$ if and only if $\phi_t((\t_{\rho}^{-1}\t A,\t_{e_{u}}))=0$ for all $u\in W_0$.
 By the formula (6) established in the proof of Theorem 3.4, we have
 $$(\t_{\rho}^{-1}\t A,\t_{e_u})=(-1)^\nu q^{-\frac{\nu}2}\sum_{w\in W_0}w(\t_{ e_u}\t\prod_{\a\in
R^+}\frac{1-q\t_\a}{1-\t_\a}).$$ Hence the condition $\phi_t((\t_{\rho}^{-1}\t A,\t_{e_{u}}))=0$ for all $u\in W_0$ is exactly the condition (8). We proved (a).

Since $C'\t_{\rho^{-1}}C=(-1)^{l(w_0)}C'\t_\rho C$, (b) follows from (a) by applying the involution $h\to h^*$. The theorem is proved.

\subsection*{4.2} The left ideal $\mathcal I$ defined in 2.1 is in the kernel of the $H$-module homomorphism $\psi: H\to H_tC', h\to hC'$. Thus $\psi$ induces an $H$-module homomorphism $ N\to H_tC'$, denoted again by $\psi$. Denote the image of $C_w$ in $H/\mathcal I$ by the same notation $C_w$. Then $C_w, \ w\in Y_0$, form a basis of $N=H/\mathcal I$.

For each two-sided cell $c$, let $N^{\le c}$ be the submodule of $N$ spanned   by all $C_w, \ w\in Y_0$ and $w\le_{L} u$ for some $u\in c\cap Y_0$. Also let $N^{< c}$ be the submodule of $N$ spanned by all $C_w, \ w\in Y_0-Y_0\cap c$ and $w\le_{L} u$ for some $u\in c\cap Y_0$. Then $\psi(N^{<c})$ and $\psi(N^{\le c})$ are submodules of $H_tC'$. This gives some natural submodules of $H_tC'$. Theorem 4.1 asserts that $\psi(N^{\le c_0})=H_tC\theta_\rho C'$ is either 0 or an irreducible submodule of $H_tC'$. In general, $\psi(N^{\le c})/\psi(N^{< c})$ maybe 0 or reducible, but it is not clear whether this module is  semisimple if it is not 0. (We refer to [KL1] for the definition of preorder $\le_L$.)

When $k$ is the field of complex numbers, by [X4, Theorem 7.8], $\psi(N^{\le c_0})=H_tC\t_\rho C'$ is  0 if and only if the set
$$\mathbf g_{t,q}=\{ X\in \text{Lie}(G)\,|\, \text{Ad}(t)(X)=qX\}$$ contains nonzero semisimple elements.

 For a simple $H$-module $L$, there exists a unique two-sided cell $c$ of $W$ such that  $C_wL\ne 0$ for some $w\in c$ and $C_uL=0$ for any $u\in W-c$ with $u\le_{LR} w$ . The two-sided cell $c$ is denoted by $c_L$ and is called the two-sided cell attached to $L$. (We refer to [KL1] for the definition of preorder $\le_{LR}$.)

 We have seen that for a simple $H$-module $L$, $c_L=c_0$ if and only if $CL\ne 0$. Theorem 4.1 gives a computable (in principle) realization of irreducible $H$-modules with attached two-sided cell $c_0$.

\bibliographystyle{unsrt}

\end{document}